\documentclass{article}
\usepackage[leqno,intlimits]{amsmath}
\usepackage{amssymb}
\usepackage{amsthm}
\usepackage{hhline}
\usepackage{hyperref} 

\newtheorem{tm}{Theorem}[section]
\newtheorem{lm}[tm]{Lemma}
\newtheorem{pr}[tm]{Proposition}
\newtheorem{cor}[tm]{Corollary}
\newtheorem{rem}[tm]{Remark}

\newcommand*{\un}[1]{\underline{#1}}
\newcommand*{\Zb}{\mathbb Z}
\newcommand*{\Rb}{\mathbb R}
\newcommand*{\om}{\omega}
\newcommand*{\ze}{\zeta}
\newcommand*{\de}{\delta}
\newcommand*{\omin}{\om^{\text{min}}}
\newcommand*{\omax}{\om^{\text{max}}}
\newcommand*{\ba}{\begin{aligned}}
\newcommand*{\ea}{\end{aligned}}
\newcommand*{\be}{\begin{equation}}
\newcommand*{\ee}{\end{equation}}
\newcommand*{\te}{\theta}
\newcommand*{\e}[1]{\text{\rm e}^{#1}}
\newcommand*{\vp}{\varphi}
\newcommand*{\vr}{\varrho}
\newcommand*{\Ev}{{\bf E}}
\newcommand*{\Pv}{{\bf P}}
\newcommand*{\Vv}{{\text{\bf Var}}}
\newcommand*{\Cov}{{\text{\bf Cov}}}
\newcommand*{\Ic}{\mathcal I}
\newcommand*{\lc}{\lceil}
\newcommand*{\rc}{\rceil}
\newcommand*{\lf}{\lfloor}
\newcommand*{\rf}{\rfloor}
\newcommand*{\ta}{{\lf Vt\rf}}

\newcommand*{\tr}{\widetilde{r}}
\newcommand*{\di}{\,\text{\rm d}}
\newcommand*{\Ff}{\,\mathcal F}
\newcommand*{\al}{\alpha}
\newcommand*{\tom}{\widetilde{\omega}}
\newcommand*{\ga}{\gamma}
\newcommand*{\wt}{\widetilde}
\newcommand*{\wih}{\widehat}
\newcommand*{\wh}{\widetilde h}
\newcommand*{\hop}{\bigskip\noindent}
\newcommand*{\up}{\uparrow}
\newcommand*{\dn}{\downarrow}

\numberwithin{equation}{section}

\begin{document}
\title{Exact connections between current fluctuations and the second class particle in a class of deposition models}
\author{M.\ Bal\'azs\thanks{M. Bal\'azs was partially supported by the Hungarian Scientific Research Fund (OTKA) grants TS49835, T037685 and National Science Foundation grant DMS-0503650.}\phantom{a}\,and T.\ Sepp\"al\"ainen\thanks{T.\ Sepp\"al\"ainen was partially supported by National Science Foundation grant DMS-0402231.}\\
University of Wisconsin-Madison}
\maketitle
\begin{abstract}
We consider a large class of nearest neighbor attractive stochastic interacting systems that includes the asymmetric simple exclusion, zero range, bricklayers' and the symmetric K-exclusion processes. We provide exact formulas that connect particle flux (or surface growth) fluctuations to the two-point function of the process and to the motion of the second class particle. Such connections have only been available for simple exclusion where they were of great use in particle current fluctuation investigations.
\end{abstract}

\hop

\section{Introduction}

Serious research has recently been conducted on particle flux fluctuations in the asymmetric simple exclusion process (ASEP). This process is one of the simplest interacting particle systems, where particles are located at sites of $\Zb$, and each of them attempts to jump a unit step to the right after an independent exponential time with rate $0\le p\le1$, and to the left with rate $q=1-p$. If the destination site is occupied by another particle at the time of the attempt then the jump is suppressed. The translation-invariant extremal stationary measures for the ASEP are the Bernoulli distributions with density $0\le\vr\le1$.

In the Eulerian scaling, when time and space are rescaled with the same parameter, the hydrodynamic limit of the process leads to the inviscid Burgers equation, see e.g.\ Rezakhanlou \cite{hl} or Spohn \cite{spohnhdl}. Burgers equation possesses a characteristic speed $V^\vr$ for each $\vr$ value of the local Bernoulli equilibrium.

Ferrari and Fontes \cite{se} considered the time-integrated particle flux (or the net number of particles) that cross by time $t$ the path of an observer who moves with a constant speed $V$. They proved that this random quantity is asymptotically normal and its variance scales with the first power of time when the observer's speed $V$ differs from the characteristic speed $V^\vr$. Fluctuations in this case basically come from the initial Bernoulli distribution. 

However, the correct scaling for integrated flux variance becomes $t^{2/3}$ when the observer moves with the characteristic speed $V=V^\vr$. This is the case when the genuine dynamical fluctuations of the process become visible. Many interesting results have been discovered in this direction, we refer to Pr\"ahofer and Spohn \cite{spohn}, Ferrari and Spohn \cite{ferspohn}, Bal\'azs, Cator and Sepp\"al\"ainen \cite{third}, Quastel and Valk\'o \cite{quava}, Bal\'azs and Sepp\"al\"ainen \cite{se2/3}.

While the above results were true in simple exclusion, many other stochastic interacting systems are expected to have similar scaling properties. An exact connection between the flux variance and second class particles, that is derived with the help of the two-point function, proved to be essential in all of the above papers. In the present note we build on the first part of Bal\'azs \cite{fluct} to derive exact formulas that provide the same type of connection in a wide class of nearest neighbor attractive stochastic interacting systems. In \cite{fluct} only asymptotic versions of these formulas were derived and used. The class we consider follows the basic ideas of Cocozza-Thi\-vent \cite{coco}, and it includes the asymmetric simple exclusion, zero range, bricklayers' processes, and the symmetric K-exclusion processes. Investigations on particle flux fluctuations will hopefully progress beyond ASEP to involve these more general models, and the present findings can serve as building blocks for such future arguments.

{\bf Notation.} Variables $i$, $j$, $k$ and $n$ will refer to sites of $\Zb$, while $x$, $y$ and $z$ will be integers in the single-site state space $I$ to be defined below.

\subsection{The model}

The class of models described here is a generalization of the so-called misanthrope process. We use a surface growth interpretation, but many members of this class can be understood in terms of particles jumping on the one dimensional lattice. For $-\infty\le\omin\le0$ and $1\le\omax\le\infty$ (possibly infinite valued) integers, we define the single-site state space
\[
I:\,=\left\{z\in\Zb\,:\,\omin-1<z<\omax+1\right\}
\]
and the configuration space
\[
\Omega=\left\{\un\om=(\om_i)_{i\in\Zb}\ :\ \om_i\in I\right\}=I^{\Zb}.
\]
For each pair of neighboring sites $i$ and $i+1$ of $\Zb$, we can imagine a column built of bricks above the edge $(i,\,i+1)$. The height of this column is denoted by $h_i$. A state configuration $\un{\om}\in\Omega$ has components $\om_i=h_{i-1}-h_i\,\in I$, being the negative discrete gradients of the height of the ``wall''. The evolution is described by jump processes. A brick can be added:
\[
\left.
\ba
\left(\om_i,\,\om_{i+1}\right)&\longrightarrow\left(\om_i-1,\,\om_{i+1}+1\right)\\
h_i&\longrightarrow h_i+1
\ea
\right\}
\text{with rate}\ p(\om_i,\,\om_{i+1}),
\]
or removed:
\[
\left.
\ba
\left(\om_i,\,\om_{i+1}\right)&\longrightarrow\left(\om_i+1,\,\om_{i+1}-1\right)\\
h_i&\longrightarrow h_i-1
\ea
\right\}\text{with rate}\ q(\om_i,\,\om_{i+1}).
\]
Conditionally on $\un\om(t)$, these moves are independent. We impose the following assumptions on the rates:
\begin{itemize}
\item The rates must satisfy
\be
p(\omin,\,\cdot\,)\equiv p(\,\cdot\,,\,\omax)\equiv q(\omax,\,\cdot\,)\equiv q(\,\cdot\,,\,\omin)\equiv0\label{eq:ominmax}
\ee
whenever either $\omin$ or $\omax$ is finite. We assume that either $p$ and $q$ are non-zero in all other cases, or one of them is the identically zero function (totally asymmetric case).
\item The dynamics will have a smoothening effect when we assume monotonicity in the following way:
\be
\ba
p(z+1,\,y)&\ge p(z,\,y),\qquad&p(y,\,z+1)&\le p(y,\,z)\\
q(z+1,\,y)&\le q(z,\,y),\qquad&q(y,\,z+1)&\ge q(y,\,z)
\ea\label{eq:mon}
\ee
for $y,\,z,\,z+1\in I$. This property has the natural interpretation that the higher neighbors a column has, the faster it grows and the slower it gets a brick removed. Our model is hence {\sl attractive}.
\item We are going to use the product property of the model's translation-invariant stationary measure. For this reason, similarly to Cocozza-Thivent \cite{coco}, we need two assumptions:
\begin{itemize}
\item For any $x,\,y,\,z\in I$
\be
\ba
&p(x,\,y)+p(y,\,z)+p(z,\,x)\!\!\!\!\!&&\\
+\,&q(x,\,y)+q(y,\,z)+q(z,\,x)\!\!\!\!\!&=&\,p(x,\,z)+p(z,\,y)+p(y,\,x)\\
&&+&\,q(x,\,z)+q(z,\,y)+q(y,\,x).
\ea\label{eq:stacifelt}
\ee
\item There are symmetric functions $s_p$ and $s_q$, and a common function $f$, such that $f(\omin)=0$ whenever $\omin$ is finite, and for any $y,\,z\in I$
\be
p(y,\,z)=s_p(y,\,z+1)\cdot f(y)\qquad\text{and}\qquad q(y,\,z)=s_q(y+1,\,z)\cdot f(z).\label{eq:symm}
\ee
Condition \eqref{eq:mon} implies that $f$ is non-decreasing on $I$.
\end{itemize}
\item In order to properly construct the dynamics, restrictive growth conditions might be necessary on the rates $p$ and $q$ in case of an unbounded single-site state space $I$. We comment on this below. In particular, we assume that all moments of the growth rates are finite w.r.t.\ the distributions $\mu_\te$ introduced in Section \ref{sc:gibbs}.
\end{itemize}
At time $t$, the interface mentioned above is described by $\un{\om}(t)$. Let $\vp\,:\,\Omega\to\Rb$ be a finite cylinder function i.e.\ $\vp$ depends on a finite number of $\om_i$ values. The growth of this interface is a Markov process, with the formal infinitesimal generator $L$:
\be
\ba
(L\vp)(\un\om)&=\sum_{i\in\Zb}p(\om_i,\,\om_{i+1})\cdot\left[\vp(\dots,\,\om_i-1,\,\om_{i+1}+1,\,\dots)-\vp(\un\om)\right]\\
&+\sum_{i\in\Zb}q(\om_i,\,\om_{i+1})\cdot\left[\vp(\dots,\,\om_i+1,\,\om_{i+1}-1,\,\dots)-\vp(\un\om)\right].
\ea\label{eq:gen}
\ee
The construction of dynamics is available in the following situations. Several models with bounded rates are well understood and can be handled via the Hille-Yosida Theorem, see Liggett \cite{ips}. When the rates $p$ and $q$ grow at most linearly fast as functions of the local $\om$ values, then methods initiated by Liggett and Andjel lead to the construction of some zero range type systems (Andjel \cite{and}, Liggett \cite{lize},  Booth and Quant \cite{lna, lorna}). The totally asymmetric zero range and bricklayers' processes with at most exponentially growing rates are constructed in Bal\'azs, Rassoul-Agha, Sepp\"al\"ainen and Sethuraman \cite{exists}. See the definition of zero range and bricklayers' processes below.

We assume that the existence of dynamics can be established on a set of tempered configurations $\wt\Omega$ (i.e.\ configurations obeying some restrictive growth conditions), and we have the usual properties of the semigroup and the generator acting on nice functions on this set. We also assume that $\wt\Omega$ is of full measure w.r.t.\ the stationary measures defined in Section \ref{sc:gibbs}. Questions of existence of dynamics are not considered in the present paper.

\subsection{Examples}

There are three essentially different cases of these models. All of them are of nearest neighbor type.

\begin{enumerate}
\item {\bf Generalized exclusion processes} are described by our models in case both $\omin$ and $\omax$ are finite.
\begin{itemize}
\item {\bf The asymmetric simple exclusion process} introduced by F.\ Spitzer \cite{spi} is described this way by $\omin=0,\ \omax=1$, $f(z)={\bf1}\{z=1\}$,
\[
s_p(y,\,z)=p\cdot{\bf1}\{y=z=1\}\qquad\text{and}\qquad s_q(y,\,z)=q\cdot{\bf1}\{y=z=1\},
\]
where $p$ and $q$ are non-negative reals adding up to 1 (see \eqref{eq:symm}). In this case
\[
p(y,\,z)=p\cdot{\bf1}\{y=1,\,z=0\}\qquad\text{and}\qquad q(y,\,z)=q\cdot{\bf1}\{y=0,\,z=1\}.
\]
Here $\om_i\in\{0,\,1\}$ is the occupation number for site $i$, $p(\om_i,\,\om_{i+1})$ is the rate for a particle to jump from site $i$ to $i+1$, and $q(\om_i,\,\om_{i+1})$ is the rate for a particle to jump from site $i+1$ to $i$. These rates have values $p$ and $q$, respectively, whenever there is a particle to perform the above jumps, and there is no particle on the terminal site of the jumps. Conditions \eqref{eq:mon} and \eqref{eq:stacifelt} are also satisfied by these rates.
\item {\bf A particle-antiparticle exclusion process} is also shown to de\-mon\-stra\-te the generality of the frame described above. Let $p+q=1$, $\omin=-1,\ \omax=1$. Fix $f(-1)=0$, $f(0)=c$ ({\sl creation}), $f(1)=a$ ({\sl annihilation}) positive rates with $c\le a/2$,
\[
\ba
s_p(0,\,1)=s_p(1,\,0)&=p,&\ s_p(0,\,0)&=\frac{pa}{2c},&\ s_p(1,\,1)&=\frac{p}{2},\\
s_q(0,\,1)=s_q(1,\,0)&=q,&\ s_q(0,\,0)&=\frac{qa}{2c},&\ s_q(1,\,1)&=\frac{q}{2}
\ea
\]
and $s_p,\ s_q$ zero in all other cases. These result in rates
\[
\ba
p(0,\,0)&=pc,&\ p(0,\,-1)=p(1,\,0)&=\frac{pa}{2},&\ p(1,\,-1)&=pa,\\
q(0,\,0)&=qc,&\ q(-1,\,0)=q(0,\,1)&=\frac{qa}{2},&\ q(-1,\,1)&=qa
\ea
\]
and zero in all other cases. If $\om_i$ is the number of particles at site $i$, with $\om_i=-1$ meaning the presence of an antiparticle, then this model describes an asymmetric exclusion process of particles and antiparticles with annihilation and particle-antiparticle pair creation. These rates also satisfy our conditions.
\end{itemize}
Other generalizations are possible allowing a bounded number of particles (or antiparticles) per site.
\item {\bf Generalized misanthrope processes} are obtained by choosing $\omin>-\infty,\ \omax=\infty$.
\begin{itemize}
\item {\bf The zero range process} is included by $p+q=1$, $\omin=0,\ \omax=\infty$, an arbitrary nondecreasing function $f\,:\,\Zb^+\to\Rb^+$ such that $f(0)=0$,
\[
\ba
s_p(y,\,z)&\equiv p\qquad&&\text{and}&\qquad s_q(y,\,z)&\equiv q,\\
p(y,\,z)&=p\cdot f(y)\qquad&&\text{and}&\qquad q(y,\,z)&=q\cdot f(z).
\ea
\]
Again, $\om_i$ represents the number of particles at site $i$. Depending on this number, a particle jumps from $i$ to the right with rate $p\cdot f(\om_i)$, and to the left with rate $q\cdot f(\om_i)$. These rates trivially satisfy conditions \eqref{eq:mon} and \eqref{eq:stacifelt}.
\end{itemize}
\item {\bf General deposition processes} are the type of these models where $\omin=-\infty$ and $\omax=\infty$. In this case, the height difference between columns next to each other can be arbitrary in $\Zb$. Hence the presence of antiparticles cannot be avoided when trying to give a particle representation of the process.
\begin{itemize}
\item {\bf Bricklayers' models.} Let $f\,:\,\Zb\to\Rb^+$ be non-decreasing, also having the property
\[
f(z)\cdot f(1-z)=1\qquad\text{for all }z\in\Zb.
\]
The values of $f$ for positive $z$'s thus determine the values for non-positive $z$'s. Define, with non-negative numbers $p+q=1$,
\[
s_p(y,\,z)=p+\frac{p}{f(y)f(z)}\qquad\text{and}\qquad s_q(y,\,z)=q+\frac{q}{f(y)f(z)},
\]
which result in
\[
p(y,\,z)=pf(y)+pf(-z)\qquad\text{and}\qquad q(y,\,z)=qf(-y)+qf(z).
\]
This process can be represented by bricklayers standing at each site $i$, laying a brick on the column on their left with rate $pf(-\om_i)$ and laying a brick to their right with rate $pf(\om_i)$. They also remove a brick from their left with rate $qf(\om_i)$ and from their right with rate $qf(-\om_i)$. This interpretation gives reason to call these models bricklayers' model. Conditions \eqref{eq:mon} and \eqref{eq:stacifelt} hold for the rates.
\end{itemize}
\item {\bf Symmetric processes} are defined to have the identity $p(y,\,z)=q(z,\,y)$. In this case \eqref{eq:stacifelt} holds automatically, we only need to take care of \eqref{eq:mon} and \eqref{eq:symm}.
\begin{itemize}
\item {\bf The symmetric K-exclusion process} is obtained if we set $\omin=0$, $\omax=K$, $f(z)={\bf1}\{z>0\}$,
\[
s_p(y,\,z)=s_q(y,\,z)={\bf1}\{z,\,y\le K\}.
\]
These result in
\[
p(y,\,z)=q(z,\,y)={\bf1}\{y>0,\ z<K\}.
\]
This process thus also has a product stationary distribution, as described below.
\end{itemize}
\end{enumerate}

\subsection{Translation invariant stationary product distributions}\label{sc:gibbs}

We now present some translation invariant stationary distributions for these processes. For many cases it has been proved that these are the only extremal translation-invariant stationary distributions. Following some ideas in Cocozza-Thivent \cite{coco}, we first consider the non-decreasing function $f$ of \eqref{eq:symm}.
 For $I\ni z>0$ we define 
\[
f(z)!:\,=\prod_{y=1}^zf(y),
\]
while for $I\ni z<0$ let
\[
f(z)!:\,=\frac{1}{\prod\limits_{y=z+1}^0f(y)},
\]
finally $f(0)!:\,=1$. Then we have 
\[
f(z)!\cdot f(z+1)=f(z+1)!
\]
for all $z\in I$. Let
\[
\bar\te:\,=\left\{\begin{array}{ll}\log\left(\liminf\limits_{z\to\infty}\left(f(z)!\right)^{1/z}\right)=\lim\limits_{z\to\infty}\log(f(z))\ \ &,\ \text{if}\ \omax=\infty\\\infty\ \ &,\ \text{else}\end{array}\right.
\]
and
\[
\un\te:\,=\left\{\begin{array}{ll}\log\left(\limsup\limits_{z\to\infty}\left(f(-z)!\right)^{1/z}\right)=\lim\limits_{z\to\infty}\log(f(-z))\ \ &,\ \text{if}\ \omin=-\infty\\-\infty\ \ &,\ \text{else}.\end{array}\right.
\]
By monotonicity of $f$, we have $\bar\te\ge\un\te$. We assume $\bar\te>\un\te$. With a generic real parameter $\te\in\left(\un\te,\,\bar\te\right)$, we define the state sum
\[
Z(\te):\,=\sum_{z\in I}\frac{\e{\te z}}{f(z)!}<\infty.
\]
Let the product-distribution $\un\mu_\te$ have marginals
\[
\mu_\te(z)=\un\mu_\te\left\{\un\om\,:\,\om_i=z\right\}:\,=\left\{
\ba
&\frac{1}{Z(\te)}\cdot\frac{\e{\te z}}{f(z)!}&&\text{if }z\in I,\\
&0&&\text{if }z\notin I.
\ea
\right.
\]
\begin{lm}\label{lm:stati}
The product distribution $\un\mu_\te$ is stationary for the process generated by \eqref{eq:gen}.
\end{lm}
We prove this lemma in Section \ref{sc:eq}.

\section{Results}

Introduce $[x]:\,=\lf x\rf$ (floor) when $x\ge0$, and $[x]:\,=\lc x\rc$ (ceil) when $x<0$. We start our process in the above described translation-invariant stationary distribution $\un\mu_\te$. The quantities of main interest are the following. For a fixed speed value $V$ we define
\[
J^{(V)}(t):=h_{[Vt]}(t)-h_0(0),
\]
the height of the column at site $[Vt]$ at time $t$, relative to the initial height of the column at the origin. For $V=0$ we write
\be
J(t)=J^{(0)}(t):\,=h_0(t)-h_0(0).\label{eq:j}
\ee
In particle notations of the models, $J^{(V)}(t)$ is the time-integrated particle flux, i.e.\ the net number of particles jumping through the moving window positioned at $Vs$, when $s$ runs in the time interval $[0,\,t]$. Often we will make the choice $h_0(0)=0$.
\begin{tm}\label{tm:covar}
For any $V\in\Rb$ and $t>0$,
\be
\Vv(J^{(V)}(t))=\sum_{n=-\infty}^\infty|[Vt]-n|\cdot\Cov(\om_n(t),\,\om_0(0)).\label{eq:main}
\ee
Moreover, we also have
\be
\sum_{n=-\infty}^\infty n\cdot\Cov(\om_n(t),\,\om_0(0))=t\cdot\Cov(p(\om_0,\,\om_1)-q(\om_0,\,\om_1),\,(\om_0+\om_1)).\label{eq:e2nd}
\ee
\end{tm}

This theorem is proved in Sections \ref{sc:eq}, \ref{sc:vg}, and \ref{sc:nv}. Our primary objects of investigation are attractive systems. There is, however, only a minor point where our arguments use attractivity, and that is at the very end of Section \ref{sc:cov2nd} where positivity of space-time correlations is used. This suggests that the above theorem should still hold in this class without attractivity, but we have not investigated this issue. However, the following results do genuinly require attractiveness.

We introduce the notation $\un\de_i\in\Omega$, a configuration being one at site $i$ and zero at all other sites. Let $\un\om\in\wt\Omega$ be such that $\om_0<\omax$, and $\un\om^+:\,=\un\om+\un\de_0$. We say in this situation that we have a single second class particle between $\un\om^+$ and $\un\om$ at site 0. Section \ref{sc:bas} describes the basic coupling of two processes. With the above initial data this coupling conserves the single second class particle for all times $t>0$:
\begin{equation}
\un\om^+(t)=\un\om(t)+\un\de_{Q(t)}.\label{eq:delt}
\end{equation}
The quantity $Q(t)$ is the position of the second class particle at time $t$, which performs a nearest neighbor walk, influenced by the ambient process $\un\om(\cdot)$.

We intend to consider the initial state $\un\om$ in distribution $\un\mu_\te$, together with a second class particle started from the origin. When $\omax<\infty$, this leads to a positive probability of $\om_0=\omax$, in which case we cannot define our initial state $\un\om^+$ in $\Omega$. We therefore introduce the indicator
\[
\wih{\bf1}\{\cdot\}:\,={\bf1}\{\cdot\}\cdot{\bf1}\{\om_0<\omax\}.
\]
With this notation it makes sense to have the $\wih{\bf1}$-indicator of a second class particle-related event under an expectation.
\begin{tm}\label{tm:cov2nd}
We use the convention that the empty sum equals zero. Let
\[
g(z):\,=z-\sum_{y\in I}y\,\mu_\te(y).
\]
For any $n\in\Zb$ and $t\ge0$, we have
\be
\Cov(\om_n(t),\,\om_0(0))=\Ev\Bigl(\wih{\bf1}\{Q(t)=n\}\cdot\sum_{z=\om_0+1}^{\omax}g(z)\frac{\mu_\te(z)}{\mu_\te(\om_0)}\Bigr).\label{eq:cov2nd}
\ee
Moreover, the sum in the last display is non-negative for any \(\om_0\).
\end{tm}
A short calculation shows that the mean of the sum on the right hand-side is $\Vv(\om_0)$ (the variance w.r.t\ the distribution $\mu_\te$). Denote by \(\un{\wih\mu}_\te\) the product measure of marginals $\un\mu_\te$ for all sites, except for the origin where
\be
\un{\wih\mu}_\te\{\un\om\,:\,\om_0=y\}:\,=\frac{1}{\Vv(\om_0)}\sum_{z=y+1}^{\omax}g(z)\mu_\te(z)\qquad(y\in I).\label{eq:modis}
\ee
We write \(\wih\Pv\) and, correspondingly, \(\wih\Ev,\ \wih\Cov\) for probabilities of a process that is started in distribution \(\un{\wih\mu}_\te\). With this notation Theorem \ref{tm:cov2nd} rewrites as
\[
\Cov(\om_n(t),\,\om_0(0))=\Vv(\om_0)\cdot\wih\Pv\{Q(t)=n\}.
\]
\begin{cor}\label{cr:poscorr}
For any $n\in\Zb$ and $t\in\Rb$, the variables $\om_0(0)$ and $\om_n(t)$ are non-negatively correlated.
\end{cor}
\begin{cor} Equation \eqref{eq:main} can be rewritten as
\be
\Vv(J^{(V)}(t))=\Vv(\om_0)\cdot\wih\Ev(|Q(t)-[Vt]|).\label{eq:modmain}
\ee
\end{cor}
Define now the quantity \(\vr=\vr_\te=\Ev_\te(\om_0)\) which is a function of the parameter \(\te\). In particle systems this is simply the density of particles. Due to the definition of \(\mu_\te\), we have \(\di\vr/\di\te=\Vv(\om_0)>0\), which shows a one-to-one correspondence between \(\vr\) and \(\te\). Define also the hydrodynamic flux \(\mathcal H(\vr)=\Ev_\te[p(\om_0,\,\om_1)-q(\om_0,\,\om_1)]\) with the value of \(\te\) for which \(\vr=\Ev(\om_0)=\Ev_\te(\om_0)\) holds. The Eulerian scaling of these processes in many cases leads to Burgers-type hydrodynamic equations of the form
\[
\partial_t\vr(T,\,X)+\partial_x\mathcal H(\vr(T,\,X))=0,
\]
where \(T\) and \(X\) are the rescaled time and space parameters (see e.g.\ Rezakhanlou \cite{hl} or Spohn \cite{spohnhdl}). It is well known that the characteristic speed of this PDE is \(V^\vr=\di\mathcal H(\vr)/\di\vr\) which, by our definitions, can be shown to equal \(\Cov(p(\om_0,\,\om_1)-q(\om_0,\,\om_1),\,(\om_0+\om_1))/\Vv(\om_0)\). 
\begin{cor}\label{cr:chsp} With our new quantities, Equation \eqref{eq:e2nd} takes the form
\be
\wih\Ev_\te(Q(t))=t\cdot V^{\vr_\te}.\label{eq:chsp}
\ee
\end{cor}

\begin{rem}
For the simple exclusion process $\un\mu_\te$ becomes the Bernoulli distribution with density $\vr=\e{\te}/(1+\e{\te})$, while \(\un{\wih\mu}_\te\) is the same distribution conditioned on \(\om_0=0\) (which is the only way of initially having a second class particle at the origin).
\end{rem}

Equation \eqref{eq:modmain} has appeared for the totally asymmetric simple exclusion process, with the special value $V=V^\vr$ (the characteristic speed), as Equation (1.7) in Ferrari and Fontes \cite{se}.

Corollary \ref{cr:chsp} has appeared in Pr\"ahofer and Spohn \cite{spohn} for the totally asymmetric simple exclusion process. It is, on one hand, in accordance with the general phenomenon that the second class particle follows the characteristics of the hydrodynamic equation. On the other hand, \(\un{\wih\mu}_\te\) is known \emph{not} to be stationary as seen by the second class particle in many models, which makes the validity of this finite-time formula somewhat surprising.

Formulas \eqref{eq:modmain} and \eqref{eq:chsp} have a great potential to serve as a basic building block in computing the order of flux fluctuations and diffusivity in our class of systems (see e.g.\ Quastel and Valk\'o \cite{quava} for the definition of diffusivity in simple exclusion), as happened in the case of simple exclusion (Quastel and Valk\'o \cite{quava}, Bal\'azs and Sepp\"al\"ainen \cite{se2/3}).

\section{Equilibrium and the reversed chain}\label{sc:eq}

This section contains some basic computations regarding the measure $\un\mu_\te$.
For $\omin-1<z,\,y<\omax+1$, define
\be
p^*(y,\,z):\,=p(z,\,y),\qquad q^*(y,\,z):\,=q(z,\,y)\label{eq:revr}
\ee
and
\be
\ba
(L^*\psi)(\un\om)&=\sum_{i\in\Zb}p^*(\om_i,\,\om_{i+1})\cdot\left[\psi(\dots,\,\om_i+1,\,\om_{i+1}-1,\,\dots)-\psi(\un\om)\right]\\
&+\sum_{i\in\Zb}q^*(\om_i,\,\om_{i+1})\cdot\left[\psi(\dots,\,\om_i-1,\,\om_{i+1}+1,\,\dots)-\psi(\un\om)\right]
\ea\label{eq:gen*}
\ee
on bounded cylinder functions $\psi$.
\begin{pr}
We have
\[
\Ev\left(\psi(\un\om)\cdot L\vp(\un\om)\right)=\Ev\left(\vp(\un\om)\cdot L^*\psi(\un\om)\right)
\]
on bounded cylinder functions $\psi$ and $\vp$.
\end{pr}

\begin{proof}
By definition we have the property that for all $\omin-1<z<\omax$,
\[
\frac{\mu_\te(z+1)}{\mu_\te(z)}=\frac{\e{\te}}{f(z+1)}.
\]
This, together with \eqref{eq:symm}, implies that whenever $\omin-1<y<\omax$ and $\omin<z<\omax+1$,
\[
p(y+1,\,z-1)\cdot\frac{\mu_\te(y+1)\,\mu_\te(z-1)}{\mu_\te(y)\,\mu_\te(z)}=p(z,\,y)
\]
holds, and whenever $\omin<y<\omax+1$ and $\omin-1<z<\omax$, then
\[
q(y-1,\,z+1)\cdot\frac{\mu_\te(y-1)\,\mu_\te(z+1)}{\mu_\te(y)\,\mu_\te(z)}=q(z,\,y)
\]
holds. With \eqref{eq:ominmax}, changing variables leads to
\begin{multline}
\sum_{y=\omin}^{\omax}\sum_{z=\omin}^{\omax}p(y,\,z)\cdot G(y-1,\,z+1)\mu_\te(y)\mu_\te(z)\\
\ba
&=\sum_{y=\omin+1}^{\omax}\sum_{z=\omin}^{\omax-1}p(y,\,z)\cdot G(y-1,\,z+1)\mu_\te(y)\mu_\te(z)\\
&=\sum_{y=\omin}^{\omax-1}\sum_{z=\omin+1}^{\omax}p(y+1,\,z-1)\cdot\frac{\mu_\te(y+1)\,\mu_\te(z-1)}{\mu_\te(y)\,\mu_\te(z)}\cdot G(y,\,z)\mu_\te(y)\mu_\te(z)\\
&=\sum_{y=\omin}^{\omax-1}\sum_{z=\omin+1}^{\omax}p(z,\,y)\cdot G(y,\,z)\mu_\te(y)\mu_\te(z)\\
&=\sum_{y=\omin}^{\omax}\sum_{z=\omin}^{\omax}p(z,\,y)\cdot G(y,\,z)\mu_\te(y)\mu_\te(z)
\ea\label{eq:cvp}
\end{multline}
for any function $G$ which makes the sums convergent. In a similar fashion, we get
\begin{multline}
\sum_{y=\omin}^{\omax}\sum_{z=\omin}^{\omax}q(y,\,z)\cdot G(y+1,\,z-1)\mu_\te(y)\mu_\te(z)\\
=\sum_{y=\omin}^{\omax}\sum_{z=\omin}^{\omax}q(z,\,y)\cdot G(y,\,z)\mu_\te(y)\mu_\te(z).\label{eq:cvq}
\end{multline}

Let $\psi,\ \vp$ be bounded cylinder functions, and let $\Ic\subset\Zb$ be a finite discrete interval of which the size can be divided by three, and which contains the set
\[
\left\{i\in\Zb\,:\,\psi\ \text{or}\ \vp\ \text{depends on}\ \om_i\ \text{or on}\ \om_{i-1}\right\}.
\]
Then the summation index $i$ in the definition \eqref{eq:gen} of the generator can be run on the set $\Ic$. We begin by changing variables $\om_i,\,\om_{i+1}$ as in \eqref{eq:cvp} and \eqref{eq:cvq}:
\[
\ba
\Ev\left(\psi(\un\om)\cdot L\vp(\un\om)\right)&=\Ev\sum_{i\in\Ic}p(\om_{i+1},\,\om_i)\cdot\psi(\dots,\,\om_i+1,\,\om_{i+1}-1,\,\dots)\vp(\un\om)\\
&\qquad-\Ev\sum_{i\in\Ic}p(\om_i,\,\om_{i+1})\cdot\psi(\un\om)\vp(\un\om)\\
&+\Ev\sum_{i\in\Ic}q(\om_{i+1},\,\om_i)\cdot\psi(\dots,\,\om_i-1,\,\om_{i+1}+1,\,\dots)\vp(\un\om)\\
&\qquad-\Ev\sum_{i\in\Ic}q(\om_i,\,\om_{i+1})\cdot\psi(\un\om)\vp(\un\om).
\ea
\]
Since $|\mathcal I|$ is divisible by three, \eqref{eq:stacifelt} implies
\[
\sum_{i\in\mathcal I}\bigl[p(\om_i,\,\om_{i+1})+q(\om_i,\,\om_{i+1})\bigr]=\sum_{i\in\mathcal I}\bigl[p(\om_{i+1},\,\om_i)+q(\om_{i+1},\,\om_i)\bigr].
\]
We thus conclude
\begin{multline*}
\Ev\left(\psi(\un\om)\cdot L\vp(\un\om)\right)\\
\ba
&=\Ev\sum_{i\in\Ic}p(\om_{i+1},\,\om_i)\cdot\bigl[\psi(\dots,\,\om_i+1,\,\om_{i+1}-1,\,\dots)-\psi(\un\om)\bigr]\vp(\un\om)\\
&+\Ev\sum_{i\in\Ic}q(\om_{i+1},\,\om_i)\cdot\bigl[\psi(\dots,\,\om_i-1,\,\om_{i+1}+1,\,\dots)-\psi(\un\om)\bigr]\vp(\un\om).
\ea
\end{multline*}
Comparing this display with \eqref{eq:gen*} finishes the proof.
\end{proof}

\begin{proof}[Proof of Lemma \ref{lm:stati}]
The previous proposition with $\psi(\un\om)\equiv1$ shows that the expectation of the generator on any bounded cylinder $\vp$ is zero.
\end{proof}

\begin{cor}
Formula \eqref{eq:gen*} is the generator of the reversed process which has rates \eqref{eq:revr}.
\end{cor}
Note that the rates of the reversed process do not depend on the parameter $\te$ of the original process' equilibrium distribution.

Define the microscopic fluxes
\[
r(y,\,z):\,=p(y,\,z)-q(y,\,z),\qquad r^*(y,\,z):\,=p^*(y,\,z)-q^*(y,\,z).
\]
The sum of the rates for a column will be
\[
S(y,\,z):\,=p(y,\,z)+q(y,\,z)
.
\]
\begin{cor}\label{cr:r*s}
\[
\Ev(r^*(\om_0,\,\om_1)\cdot(\om_0-\om_1))=-\Ev(S(\om_0,\,\om_1)).
\]
\end{cor}
\begin{proof}
The proposition and Formula \eqref{eq:cvp} implies
\[
\ba
\Ev(p^*(\om_0,\,\om_1)\cdot(\om_0-\om_1))&=\Ev(p(\om_1,\,\om_0)\cdot(\om_0-\om_1))\\
&=\Ev(p(\om_0,\,\om_1)\cdot(\om_0-\om_1-2))\\
&=-\Ev(p^*(\om_0,\,\om_1)\cdot(\om_0-\om_1))-2\Ev(p(\om_0,\,\om_1)),
\ea
\]
from which $\Ev(p^*(\om_0,\,\om_1)\cdot(\om_0-\om_1))=-\Ev(p(\om_0,\,\om_1))$. Similarly,
\[
\ba
\Ev(q^*(\om_0,\,\om_1)\cdot(\om_0-\om_1))&=\Ev(q(\om_1,\,\om_0)\cdot(\om_0-\om_1))\\
&=\Ev(q(\om_0,\,\om_1)\cdot(\om_0-\om_1+2))\\
&=-\Ev(q^*(\om_0,\,\om_1)\cdot(\om_0-\om_1))+2\Ev(q(\om_0,\,\om_1)),
\ea
\]
thus $\Ev(q^*(\om_0,\,\om_1)\cdot(\om_0-\om_1))=\Ev(q(\om_0,\,\om_1))$.
\end{proof}

\section{Vertical growth}\label{sc:vg}

It is easier to first consider the variance of $J(t)$ (case $V=0$).

\subsection{Martingale tricks}

For convenience, we introduce the notation $\wt A$ for the centered random variable $A\in\text{L}^1(\un\mu_\te)$, and we further simplify notations by
\[
r(t):\,=r(\om_0(t),\,\om_1(t)),\quad r^*(t):\,=r^*(\om_0(t),\,\om_1(t)),\quad S(t):\,=S(\om_0(t),\,\om_1(t)).
\]
Recall \eqref{eq:j}.
\begin{lm}\label{lm:mart}
\[
\Vv(J(t))=t\,\Ev(S)+2\int_0^t\int_0^s\,\Ev\left(\tr(v)\,r^*(0)\right)\di v\,\di s.
\]
\end{lm}
\begin{proof}
First notice that applying the generator on $J^k(t)$, then using H\"older's inequality with the assumption that all moments of the rates $p$ and $q$ are finite implies in an inductive fashion that (the time-derivative of) any moment of $J(t)$ is finite. By definition, $\Ev(J(t)\,|\,\un\om(0))=t\,r(0)+\mathfrak{o}(t)$, hence
\be
M(t):\,=J(t)-\int_0^tr(s)\di s\label{eq:mtg}
\ee
is a martingale with $M(0)=0$. Therefore,
\begin{equation}
\Vv(J(t))=\Ev M(t)^2+2\,\Ev\Bigl(M(t)\,\int_0^t\tr(s)\di s\Bigr)+\Ev\Bigl(\Bigl(\int_0^t\tr(s)\di s\Bigr)^2\Bigr).\label{eq:rovidebb}
\end{equation}
Due to $\Ev\left(M(t)^2\,|\,\un\om(0)\right)=t\,S(0)+\mathfrak{o}(t)$, the process
\[
N(t):\,=M(t)^2-\int_0^t S(s)\di s
\]
is also a martingale with $N(0)=0$. Hence
\[
\Ev M(t)^2=t\,\Ev(S).
\]
Using the martingale property of $M$, the second term of \eqref{eq:rovidebb} can be written as
\[
2\,\int_0^t\Ev\left(M(t)\,\tr(s)\right)\di s=2\int_0^t\Ev\left(M(s)\,\tr(s)\right)\di s.
\]
Simply changing the limits of integration in the third term of \eqref{eq:rovidebb}, we have
\[
\Ev\Bigl(\Bigl(\int_0^t\tr(s)\di s\Bigr)^2\Bigr)=2\int_0^t\Ev\Bigl(\tr(s)\int_0^s\tr(u)\di u\Bigr)\di s.
\]
These calculations lead to
\begin{multline}
\Vv(J(t))=t\,\Ev(S)+2\int_0^t\Ev\Bigl(\tr(s)\Bigl(M(s)+\int_0^s\tr(u)\di u\Bigr)\Bigr)\di s\\
=t\,\Ev(S)+2\int_0^t\Ev\Bigl(\tr(s)\,J(s)\Bigr)\di s.\label{eq:hosszu}
\end{multline}
In order to handle $\Ev\left(\tr(s)\,J(s)\right)$, we introduce $J^{(s)\,*}$, the quantity corresponding to $J$ in the reversed process by
\[
J^{(s)\,*}(u):\,=J(s)-J(s-u)\ \ \ \ \ (s\ge u\ge0). 
\]
This is the number of bricks removed from the column in the reversed process started from time $s$. As in case of $J(t)$, a reversed martingale can be separated by
\[
M^{(s)\,*}(u):\,=J^{(s)\,*}(u)-\int_0^ur^*(s-v)\di v.
\]
For this reversed object,\,\ $M^{(s)\,*}(0)=0$\,\ and\,\ $\Ev\left(M^{(s)\,*}(u)\,|\,\Ff_{[t,\,\infty)}\right)=M^{(s)\,*}(s-t)$\ \ if\ \ $0\le s-t\le u$, where $\Ff$ stands for the natural filtration of the (forward) process. In view of this, 
\begin{multline*}
\Ev\left(\tr(s)\,J(s)\right)=\Ev\left[\tr(s)\,\Ev\left(J^{(s)\,*}(s)\,|\,\Ff_{[s,\,\infty)}\right)\right]=\\
=\Ev\left(\tr(s)\int_0^sr^*(s-v)\di v\right)=\int_0^s\Ev\left(\tr(v)\,r^*(0)\right)\di v,
\end{multline*}
where in the last step we used time-invariance of the measure. Using this result, we obtain the statement from (\ref{eq:hosszu}) by changing the order of integration.
\end{proof}

\subsection{Space-time correlations}

In this subsection we denote $r(\om_i,\,\om_{i+1})$ and $\tr(\om_i,\,\om_{i+1})$ by $r_i$ and $\tr_i$, respectively. We also keep the notation $r=r_0$ and $\tr=\tr_0$ from the previous section. For $k\in\Zb$, let 
\[
d_k\ \ :\ \ \Omega\to I\ \ ;\ \ d_k(\un\om)=\om_k
\]
be the $k$-th coordinate of $\Omega$. Then
\be
\left(Ld_k\right)(\un\om)=r_{k-1}-r_k\qquad\text{and}\qquad\left(L^*d_k\right)(\un\om)=-r^*_{k-1}+r^*_k,\label{eq:szomsz}
\ee
where $L^*$ is the infinitesimal generator (\ref{eq:gen*}) for the reversed process.
\begin{lm}\label{lm:elketto}
For $0<\al<1$ the functions
\begin{equation}
\vp_\al:\,=\sum_{k=1}^{\infty}\al^{k-1}d_k\qquad\qquad\qquad\psi_\al:\,=\sum_{k=0}^{\infty}\al^kd_{-k}\label{eq:fidef}
\end{equation}
are $\un\mu_\te$-a.s.\ well defined and finite, and
\begin{alignat*}{4}
\lim_{\al\to1}(L\vp_\al)(\un\om)&=-\lim_{\al\to1}(L\psi_\al)(\un\om)&=&\,\tr,\\
\lim_{\al\to1}(L^*\psi_\al)(\un\om)&=-\lim_{\al\to1}(L^*\vp_\al)(\un\om)&=&\,\wt r^*
\end{alignat*}
in $\text{L}^2(\un\mu_\te)$.
\end{lm}
\begin{proof}
The a.s.\ existence of the sums above can easily be shown by using the Borel-Cantelli lemma for the sets 
\[
A_n:\,=\left\{\un\om\ :\ |\om_n|\ge n\right\}.
\]
We show the first equation for $\vp_\al$. By (\ref{eq:szomsz}) 
\be
\left(L\vp_\al\right)(\un\om)=r_0+(\al-1)\,\sum_{k=1}^{\infty}r_k\al^{k-1}=\tr_0+(\al-1)\,\sum_{k=1}^{\infty}\tr_k\al^{k-1}.\label{eq:l2}
\ee
By independence of $\om_i$ and $\om_j$ for $i\ne j,\ \Ev(\tr_l\cdot\tr_k)=0$ if $|l-k|>1$ and
\[
\left|\Ev(\tr_l\cdot\tr_k)\right|=|\Cov(r_l,\,r_k)|\le\Vv(r_l)=\Ev(\tr_l\cdot\tr_l)=||\tr||_2^2,
\]
if $|k-l|=0$ or 1. Hence the $\text{L}^2$-norm of the second term on the right-hand side of (\ref{eq:l2}) tends to zero as $\al\to1$:
\begin{multline*}
\left|\left|(\al-1)\,\sum_{k=1}^{\infty}\tr_k\al^{k-1}\right|\right|_2^2
\le(\al-1)^2\sum_{k=1}^{\infty}||\tr_k||_2^2\al^{2k-2}+\\
+2(\al-1)^2\sum_{k=1}^{\infty}||\tr_k||_2^2\al^{2k-3}
=\frac{(\al-1)^2}{1-\al^2}||\tr||_2^2(1+2\al^{-1})\underset{\al\to1}{\longrightarrow}0.
\end{multline*}
The proof of the other three equations is similar.
\end{proof}
\begin{lm}\label{lm:fw}
Let $\Phi(0)=\Phi(\un\om(0))$ be an $\text{L}^2$-function that depends on the initial state only. Then
\[
\ba
\int_0^t\Ev(\tr(v)\,\Phi(0))\di v&=\lim_{\al\to1}\bigl[\Ev(\wt\vp_\al(t)\,\Phi(0))-\Ev(\wt\vp_\al(0)\,\Phi(0))\bigr]\\
&=\lim_{\al\to1}\bigl[\Ev(\wt\psi_\al(0)\,\Phi(0))-\Ev(\wt\psi_\al(t)\,\Phi(0))\bigr].
\ea
\]
\end{lm}
\begin{proof}
We show the first equality, and by convenience we also center $\Phi(0)$ besides the centered rates. We first make use of the $\text{L}^2$ convergence of the previous lemma:
\begin{multline*}
\left|\int_0^t\Ev(\tr(v)\,\wt\Phi(0))\di v-\lim_{\al\to1}\int_0^t\Ev(L\vp_\al(v)\,\wt\Phi(0))\di v\right|\\
\ba
&\le\lim_{\al\to1}\int_0^t\sqrt{\Ev\left([\tr(v)-L\vp_\al(v)]^2\right)\cdot\Ev(\wt\Phi(0)^2)}\di v\\
&=\lim_{\al\to1}t\sqrt{\Ev\left([\tr(0)-L\vp_\al(0)]^2\right)\cdot\Ev(\wt\Phi(0)^2)}=0.
\ea
\end{multline*}
Next we write
\[
\int_0^t\Ev(L\vp_\al(v)\,\wt\Phi(0))\di v=\int_0^t\Ev\Bigl(\sum_{k=1}^\infty\al^{k-1}Ld_k(\un\om(v))\,\wt\Phi(0)\Bigr)\di v.
\]
Due to \eqref{eq:szomsz}, a Cauchy inequality on $\Ev([r_{k-1}-r_k]\,\wt\Phi(0))$, finite moments of the rates and translation invariance, $\al^{k-1}Ld_k(\un\om(v))\,\wt\Phi(0)$ is absolute summable, ``expectable'' and integrable. Therefore
\[
\ba
\int_0^t\Ev(L\vp_\al(v)\,\wt\Phi(0))\di v&=\sum_{k=1}^\infty\al^{k-1}\int_0^t\Ev\bigl(Ld_k(\un\om(v))\,\wt\Phi(0)\bigr)\di v\\
&=\sum_{k=1}^\infty\al^{k-1}\bigl[\Ev\bigl(d_k(\un\om(t))\,\wt\Phi(0)\bigr)-\Ev\bigl(d_k(\un\om(0))\,\wt\Phi(0)\bigr)\bigr]\\
&=\Ev\bigl(\vp_\al(t)\,\wt\Phi(0)\bigr)-\Ev\bigl(\vp_\al(0)\,\wt\Phi(0)\bigr)
\ea
\]
by the integrated Kolmogorov equation on $d_k$ and another absolute integrability argument.
\end{proof}
\begin{lm}\label{lm:bw}
With $\Phi(0)$ as above,
\[
\ba
\int_0^t\Ev(\tr^*(-v)\,\Phi(0))\di v&=\lim_{\al\to1}\bigl[\Ev(\wt\psi_\al(-t)\,\Phi(0))-\Ev(\wt\psi_\al(0)\,\Phi(0))\bigr]\\
&=\lim_{\al\to1}\bigl[\Ev(\wt\vp_\al(0)\,\Phi(0))-\Ev(\wt\vp_\al(-t)\,\Phi(0))\bigr].
\ea
\]
\end{lm}
\begin{proof}
For the first line, repeat the previous proof except for the use of $L^*$ and $\psi_\al$ rather than $L$ and $\vp_\al$, and
\[
\int_0^t\Ev\bigl(L^*d_{-k}(\un\om(-v))\,\wt\Phi(0)\bigr)\di v=\Ev\bigl(d_{-k}(\un\om(-t))\,\wt\Phi(0)\bigr)-\Ev\bigl(d_{-k}(\un\om(0))\,\wt\Phi(0)\bigr).
\]
\end{proof}
Now we can compute the integrals in our expression for $\Vv(J)$.
\begin{tm}\label{tm:formula1}
\[
\Vv(J(t))=\sum_{n=-\infty}^{\infty}|n|\cdot\Ev(\tom_0(0)\,\tom_n(t)).
\]
\end{tm}
\begin{proof}
The aim is to rewrite the double integral of Lemma \ref{lm:mart}. A slight modification of Lemma \ref{lm:fw} that includes a second integral as well implies that this can be done by
\[
\Vv(J(t))=t\,\Ev(S)+2\lim_{\al\to1}\int_0^t\Ev(\wt\vp_\al(s)\,r^*(0))\di s-2t\lim_{\al\to1}\Ev(\wt\vp_\al(0)\,r^*(0)).
%
%
\]
The limits could be distributed to the difference since the second limit is finite by the centering and the product structure of $\un\mu_\te$. The integral in this display is rewritten with a time-translation as
\begin{multline*}
\int_0^t\Ev(\wt\vp_\al(s)\,\tr^*(0))\di s=\int_0^t\Ev(\wt\vp_\al(0)\,\tr^*(-s))\di s\\
=\lim_{\ga\to1}\Ev(\wt\vp_\al(0)\,\wt\psi_\ga(-t))-\lim_{\ga\to1}\Ev(\wt\vp_\al(0)\,\wt\psi_\ga(0)).
\end{multline*}
by Lemma \ref{lm:bw}. Notice again here that the second term on the right is finite. Hence with definitions \eqref{eq:fidef}, the variance of $J(t)$ can now be written as
\begin{align}
\Vv(J(t))&=t\,\Ev(S)-2\,t\lim_{\al\to1}\Ev(\wt\vp_\al(0)\,\tr^*(0))\notag\\
&\qquad+2\lim_{\al,\ga\to1}\Ev(\wt\vp_\al(0)\,\wt\psi_\ga(-t))-2\lim_{\al,\ga\to1}\Ev(\wt\vp_\al(0)\,\wt\psi_\ga(0))\notag\\
&=t\,\Ev(S)-2\,t\,\lim_{\al\to1}\Ev\left(\sum_{k=1}^{\infty}\al^{k-1}\,\wt\om_k(0)\,\tr^*(0)\right)\notag\\
&\qquad+2\,\lim_{\al,\ga\to1}\Ev\left(\sum_{k=1}^{\infty}\al^{k-1}\,\wt\om_k(0)\sum_{l=0}^{\infty}\ga^l\,\wt\om_{-l}(-t)\right)\label{eq:stcov}\\
&\qquad-2\,\lim_{\al,\ga\to1}\Ev\left(\sum_{k=1}^{\infty}\al^{k-1}\,\wt\om_k(0)\sum_{l=0}^{\infty}\al^l\,\wt\om_{-l}(0)\right).\notag
\end{align}
Using product property of the measure at time $t=0$ and the fact that $r^*$ depends only on $\om_0$ and $\om_1$, most of our expressions become simple (recall that all quantities with tilde are centered random variables). We prove in Section \ref{sc:cov2nd} that the limits in \eqref{eq:stcov} can be taken under the sum and the expectation:
\begin{align}
\Vv(J(t))&=t\,\Ev(S)-2\,t\,\Ev(\tom_1(0)\,\tr^*(0))+2\sum_{k=1}^\infty\sum_{l=0}^\infty\Ev(\tom_k(0)\,\tom_{-l}(-t))-0\notag\\
&=t\,\Ev(S)-2\,t\,\Ev(\tr^*(0)\,\tom_1(0))+2\sum_{n=1}^\infty n\,\Ev(\tom_n(t)\,\tom_0(0)).\label{eq:varjr}
\end{align}
In the last step, we used space and time translation-invariance of the measure.

We took advantage of the first equalities in both Lemma \ref{lm:fw} and \ref{lm:bw}. The second identities therein can be used in a similar way to prove
\be
\Vv(J(t))=t\,\Ev(S)+2\,t\,\Ev(\tr^*(0)\cdot\wt\om_0(0))+2\sum_{n=1}^{\infty}n\Ev(\tom_0(0)\,\tom_{-n}(t)).\label{eq:varjl}
\ee
The statement now follows from Corollary \ref{cr:r*s} by taking the average of the previous two displays.
\end{proof}
\begin{proof}[Proof of \eqref{eq:e2nd}]
We only need to subtract \eqref{eq:varjl} from \eqref{eq:varjr} above to conclude the last statement of Theorem \ref{tm:covar}. Notice that $r^*(y,\,z)=r(z,\,y)$, but the symmetric term $y+z$ makes the order of $y$ and $z$ in $r(z,\,y)$ immaterial in the covariance.
\end{proof}

\section{Non-vertical growth}\label{sc:nv}

We turn to the variance of the quantity $J^{(V)}(t)$. First we consider $V>0$ values.
\be
\ba
\Vv(J^{(V)}(t))&=\Ev\Bigl\{\Bigl(\wh_\ta(t)\Bigr)^2\Bigr\}
\\
&=\Ev\Bigl\{\Bigl(\wh_\ta(t)-\wh_\ta(0)\Bigr)^2\Bigr\}-\Ev\Bigl\{\Bigl(\wh_\ta(0)\Bigr)^2\Bigr\}\\
&\qquad+2\,\Ev\Bigl(\wh_\ta(t)\,\wh_\ta(0)\Bigr).\label{eq:mozgv}
\ea
\ee
Due to translation-invariance, the first term is $\Vv(J(t))$, computed in the previous sections. By $\om_i=h_{i-1}-h_i$, $h_0(0)=0$ and the product structure of the measure, the second term of the right-hand side of \eqref{eq:mozgv} is
\be
-\Ev\Bigl\{\Bigl(\wh_\ta(0)\Bigr)^2\Bigr\}=-\ta\cdot\Ev(\tom_0(0)^2)=-\ta\cdot\Vv(\om_0).\label{eq:2nd}
\ee
We compute the third term in the following lemma.
\begin{lm}\label{lm:thirdt}
For $V>0$,
\[
2\Ev\left(\wh_\ta(t)\,\wh_\ta(0)\right)=\sum_{n=-\infty}^\infty(|\ta-n|-|n|)\cdot\Ev(\tom_n(t)\,\tom_0(0))+\ta\cdot\Vv(\om_0).
\]
\end{lm}
\begin{proof}
Using $\om_i=h_{i-1}-h_i$ and $h_0(0)=0$ again,
\be
\Ev\left(\wh_\ta(t)\,\wh_\ta(0)\right)=-\sum_{j=1}^\ta\Ev(h_0(t)\,\tom_j(0))+\sum_{i=1}^\ta\sum_{j=1}^\ta\Ev(\tom_i(t)\,\tom_j(0)).\label{eq:moz3}
\ee
With the martingale of \eqref{eq:mtg} (for $J(t)=h_0(t)$), we have
\[
\Ev(h_0(t)\,\tom_j(0))=\int_0^t\Ev(r_0(s)\,\tom_j(0))\di s.
\]
We proceed via Lemma \ref{lm:fw}:
\[
\ba
\int_0^t\Ev(\tr_0(s)\,\tom_j(0))\di s&=\lim_{\al\to1}\Ev(\vp_\al(t)\,\tom_j(0))-\lim_{\al\to1}\Ev(\vp_\al(0)\,\tom_j(0))\\
&=\sum_{i=1}^\infty\Ev(\tom_i(t)\,\tom_j(0))-\Vv(\om_j(0)).
\ea
\]
Again we justify in Section \ref{sc:cov2nd} that the limit can be taken under the sum on the space-time covariances. Plugging this in \eqref{eq:moz3} gives
\be
\ba
\Ev\left(\wh_\ta(t)\,\wh_\ta(0)\right)&=-\sum_{i=\ta+1}^\infty\sum_{j=1}^\ta\Ev(\tom_i(t)\,\tom_j(0))+\ta\Vv(\om_0)\\
&=-\sum_{n=1}^\infty(n\land\ta)\Ev(\tom_n(t)\,\tom_0(0))+\ta\Vv(\om_0).
\ea\label{eq:ehhphi}
\ee
However, we also have, by the second identity of Lemma \ref{lm:fw},
\[
\ba
\int_0^t\Ev(\tr_0(s)\,\tom_j(0))\di s&=-\lim_{\al\to1}\Ev(\psi_\al(t)\,\tom_j(0))+\lim_{\al\to1}\Ev(\psi_\al(0)\,\tom_j(0))\\
&=-\sum_{i=-\infty}^0\Ev(\tom_i(t)\,\tom_j(0)).
\ea
\]
Use this in \eqref{eq:moz3} to obtain
\be
\ba
\Ev\left(\wh_\ta(t)\,\wh_\ta(0)\right)&=\sum_{i=-\infty}^\ta\sum_{j=1}^\ta\Ev(\tom_i(t)\,\tom_j(0))\\
&=\sum_{n=-\infty}^\ta(\ta-n^+)\Ev(\tom_n(t)\,\tom_0(0)).
\ea\label{eq:ehhpsi}
\ee
The statement now follows from taking the sum of \eqref{eq:ehhphi} and \eqref{eq:ehhpsi}, and from
\[
(\ta-n^+)\cdot{\bf1}\{n\le\ta\}-(n\land\ta)\cdot{\bf1}\{n\ge1\}=|\ta-n|-|n|\qquad(n\in\Zb).
\]
\end{proof}
\begin{proof}[Proof of Theorem \ref{tm:covar}]
The case $V=0$ and \eqref{eq:e2nd} are proved in the previous sections. For $V>0$, combine \eqref{eq:mozgv}, Theorem \ref{tm:formula1}, \eqref{eq:2nd}, and the previous lemma. A similar computation shows the same result for negative $V$'s.
\end{proof}

\section{The second class particle}

In this section we show how to couple a pair of our models, with the help of the so-called second class particles. The space-time correlations seen in Theorem \ref{tm:covar} are rewritten in terms of the motion of this particle.

\subsection{The basic coupling}\label{sc:bas}

For configurations $\un\eta$ and $\un\ze$ we say that $\un\eta\le\un\ze$, if $\eta_i\le\ze_i$ for all $i\in\Zb$. We consider two realizations of a process, namely, $\un\eta(\cdot)$ and $\un\ze(\cdot)$. We show the basic coupling which preserves
\begin{equation}
\un\eta(t)\le\un\ze(t)\label{eq:rend}
\end{equation}
if this property holds initially for $t=0$. We say that $d_i:\,=\ze_i(t)-\eta_i(t)\ge 0$ is the number of {\sl second class particles} present at site $i$ at time $t$.

The height of the column of $\un\ze$ (or $\un\eta$) between sites $i$ and $i+1$ is denoted by $g_i$ (or $h_i$, respectively). Let $g_i\up$ (or $h_i\up$) mean that the column of $\un\ze$ (or the column of $\un\eta$, respectively) between the sites $i$ and $i+1$ has grown by one brick. Similarly, $g_i\dn$ (or $h_i\dn$) means a brick-removal. Then the coupling rules are shown in Table \ref{tab:bas}. Each line of this table represents a possible move, with rate written in the first column. As an illustration, we also indicate the change in the number $d_i$ of second class particles at site $i$. These changes represent nearest neighbor walks of the second class particles, hence the total number of these particles is preserved.

\begin{table}[ht]
\[
\begin{array}{|c||c|c||c|c|}
\hline
\text{with rate}&g_i&h_i&d_i&d_{i+1}\\
\hhline{|=#=|=#=|=|}
p(\ze_i,\,\ze_{i+1})-p(\eta_i,\,\ze_{i+1})&\up&&\dn&\up\\
\hline
p(\eta_i,\,\eta_{i+1})-p(\eta_i,\,\ze_{i+1})&&\up&\up&\dn\\
\hline
p(\eta_i,\,\ze_{i+1})&\up&\up&&\\
\hhline{|=#=|=#=|=|}
q(\ze_i,\,\ze_{i+1})-q(\ze_i,\,\eta_{i+1})&\dn&&\up&\dn\\
\hline
q(\eta_i,\,\eta_{i+1})-q(\ze_i,\,\eta_{i+1})&&\dn&\dn&\up\\
\hline
q(\ze_i,\,\eta_{i+1})&\dn&\dn&&\\
\hline
\end{array}
\]
\caption{The basic coupling}\label{tab:bas}
\end{table}

This coupling coincides with the well-known basic coupling for particle systems. The rates of these steps are non-negative due to \eqref{eq:rend} and monotonicity \eqref{eq:mon} of $p$ and $q$. These rules clearly preserve property \eqref{eq:rend}, since the rate of any move decreasing $d_i$ becomes zero when $d_i=0$. Summing marginally the rates of jumps of either $\un\eta(\cdot)$ or $\un\ze(\cdot)$ shows that each process evolves according to its own rates.

\subsection{Space-time covariance and second class particles}\label{sc:cov2nd}

Recall the setting of Theorem \ref{tm:cov2nd}. We now show how that theorem is derived from Theorem \ref{tm:covar}. Recall that $\mu_\te(z)=0$ for any $z<\omin$ or $z>\omax$.
\begin{lm}
For the pair $(\un\om(t),\,Q(t))$ defined above Theorem \ref{tm:cov2nd} and for a function $F\,:\,I\to\Rb$ with $F(\omax)=0$ and with finite expectation value $\sum F(z)\,\mu_\te(z)$,
\begin{multline*}
\Ev\left(\om_n(t)\,\left[\frac{F(\om_0(0)-1)\,\mu_\te(\om_0(0)-1)}{\mu_\te(\om_0(0))}-F(\om_0(0))\right]\right)\\
=\Ev\left(\wih{\bf1}\{Q(t)=n\}\,F(\om_0(0))\right).
\end{multline*}
\end{lm}
\begin{proof}
Fix $\omin-1<z<\omax$, and take conditional expectation of \eqref{eq:delt}:
\begin{equation}
\Ev\left(\om^+_n(t)\,|\,\om_0(0)=z\right)=\Ev\left(\om_n(t)\,|\,\om_0(0)=z\right)+\Pv\left(Q(t)=n\,|\,\om_0(0)=z\right).\label{eq:lep}
\end{equation}
Initially, $\un\om^+(0)=\un\om(0)+\un\de_0$. Therefore, $\un\om^+(\cdot)$ itself is also a process with initial distribution $\un\mu_\te$, except for the origin. Hence
\[
\Ev\left(\om^+_n(t)\,|\,\om_0(0)=z\right)=\Ev\left(\om^+_n(t)\,|\,\om^+_0(0)=z+1\right)=\Ev\left(\om_n(t)\,|\,\om_0(0)=z+1\right),
\]
and \eqref{eq:lep} can be written as
\[
\Ev\left(\om_n(t)\,|\,\om_0(0)=z+1\right)-\Ev\left(\om_n(t)\,|\,\om_0(0)=z\right)=\Pv\left(Q(t)=n\,|\,\om_0(0)=z\right).
\]
We multiply both sides with $F(z)\,\mu_\te(z)$ and then add up for all $\omin-1<z<\omax$ to obtain
\begin{multline*}
\sum_{z\in I}\Ev\left(\om_n(t)\,|\,\om_0(0)=z\right)\cdot\left(F(z-1)\,\mu_\te(z-1)-F(z)\,\mu_\te(z)\right)\\
=\sum_{z\in I}\Pv\left(Q(t)=n\,|\,\om_0(0)=z\right)\cdot F(z)\,\mu_\te(z)
\end{multline*}
(recall $\mu_\te(z)=0$ for $z\notin I$ and $F(\omax)=0$). The proof is finished by $\Pv(\om_0(0)=z)=\mu_\te(z)$.
\end{proof}
\begin{proof}[Proof of Theorem \ref{tm:cov2nd}]
By the previous lemma, our goal is now to find the correct function $F$ with finite mean, for which $F(\omax)=0$ and
\[
\frac{F(z-1)\,\mu_\te(z-1)}{\mu_\te(z)}-F(z)=g(z)=z-\sum_{y\in I}y\,\mu_\te(y)
\]
hold. By inverting the operation on the left side, we find
\[
F(z):\,=\sum_{y=z+1}^{\omax}g(y)\,\frac{\mu_\te(y)}{\mu_\te(z)}.
\]
This function satisfies the conditions of the lemma, and \eqref{eq:cov2nd} is proved.

Now we turn to the proof that the sum in \eqref{eq:cov2nd} is non-negative. Notice that $\Ev(g(\om))=0$. The proof follows by both ${\bf1}\{z>y\}$ and $g(z)$ being non-decreasing in $z$, and hence
\[
0\le\Cov({\bf1}\{\om>y\},\,g(\om))
=\sum_{z=y+1}^{\omax}g(z)\,\mu_\te(z).
\]
\end{proof}

We finally justify taking the limits under the summations in \eqref{eq:stcov} (and at some points in the proof of Lemma \ref{lm:thirdt}). For any $\al<1$ and $\ga<1$,
\[
\sum_{k=1}^{\infty}\al^{k-1}\sum_{l=0}^{\infty}\ga^l\,\Ev(|\wt\om_k(0)|\cdot|\wt\om_{-l}(-t)|)\le\sum_{k=1}^{\infty}\al^{k-1}\sum_{l=0}^{\infty}\ga^l\,\Vv(\om_0)<\infty,
\]
and so we can reorder summations and expectations in (half of) \eqref{eq:stcov} as
\[
\sum_{k=1}^{\infty}\sum_{l=0}^{\infty}\al^{k-1}\ga^l\,\Cov(\wt\om_k(0),\,\wt\om_{-l}(-t)).
\]
Now the limits in $\al$ and $\ga$ can be brought under the double sum by Corollary \ref{cr:poscorr} and Monotone Convergence.

\section*{Acknowledgments}

M.\ Bal\'azs wishes to thank B\'alint T\'oth for initiating the study of bricklayers' models, for giving some of the basic ideas of this paper, and for helping him in many questions. We also wish to thank Christophe Bahadoran for fruitful conversations on the subject.

\bibliography{refsmarton}
\bibliographystyle{plain}

\bigskip
{\sc M.\ Bal\'azs, Department of Stochastics, Budapest University of Technology and Economics,} 1 Egry J\'ozsef u., 1111 Budapest, Hungary.\\
\indent
{\it E-mail address:} {\tt balazs@math.wisc.edu}

\bigskip
{\sc T.\ Sepp\"al\"ainen, Mathematics Department, University of Wis\-con\-sin-Madison,} Van Vleck Hall, 480 Lincoln Dr, Madison WI 53706-1388, USA.\\
\indent
{\it E-mail address:} {\tt seppalai@math.wisc.edu}

\end{document}